\documentclass{article}
\usepackage{amssymb}
\usepackage{amsfonts}
\usepackage{amsmath}

\input{tcilatex}
\begin{document}

\title{On the compactness of embeddings for a class of weighted
Orlicz-Sobolev sequence spaces}
\author{Pierre-A. Vuillermot$^{\ast ,\ast \ast }$ \\
Center for Matematical Studies, CEMS.UL, Faculdade de Ci\^{e}ncias\\
Universidade de Lisboa, 1749-016 Lisboa, Portugal$^{\ast }$\\
Universit\'{e} de Lorraine, CNRS, IECL, F-54000 Nancy, France$^{\ast \ast }$%
\\
pavuillermot@fc.ul.pt\\
ORCID Number: 0000-0003-3527-9163}
\date{}
\maketitle

\begin{abstract}
In this article we introduce a new scale of weighted Orlicz-Sobolev sequence
spaces generated by a class of suitable Orlicz functions and prove various
continuity and compactness criteria for them. In a nutshell, continuity is a
consequence of pointwise comparison between Orlicz functions while
compactness follows from the combination of the existence of a Schauder
basis in the spaces under consideration with a condition on the generating
Orlicz functions regarding their local behavior in a small neighborhood of
the origin. We illustrate our results by means of several concrete examples
and also mention some open questions along the way.

\ \ 



\ \textbf{\ \ \ \ Abbreviated title}: Compact embeddings

\ \ 
\end{abstract}

\section{Introduction and outline}

Compact embeddings for Sobolev spaces constitute an essential tool regarding
the analysis of initial-and boundary-value problems for ordinary and partial
differential equations, among other topics (see, e.g., \cite{edmundsevans}).
Extensions of such embeddings to Orlicz-Sobolev spaces defined on regions of
Euclidean space with respect to Lebesgue measure such as those described in
Theorem 8.43 of \cite{adamsfournier} or in \cite{vuillermotter} become
necessary in certain contexts as soon as exponentially growing
nonlinearities appear (see, e.g., \cite{vuillermotquarto} and some of the
references therein). To the best of our knowledge, however, no such
compactness criteria exist for Orlicz-Sobolev \textit{sequence} spaces, that
is, for Orlicz-Sobolev spaces defined with respect to discrete or counting
measures, apart from weak sequential compactness results stated in Section
10.3 of Chapter 10 in \cite{raoren2} which are relative to Orlicz sequence
spaces alone. That is the reason why we introduce here a scale of weighted
Orlicz-Sobolev sequence spaces and prove several continuous and compact
embeddings for them. The new scale generalizes the Banach spaces $h_{\mathbb{%
C},w}^{k,s}$ that were introduced in \cite{vuillermot}, consisting of all
complex sequences $\mathsf{p}=\left( p_{\mathsf{m}}\right) _{\mathsf{m}\in 
\mathbb{Z}}$ satisfying%
\begin{equation}
\sum_{\mathsf{m}\in \mathbb{Z}}w_{\mathsf{m}}\left( 1+\left\vert \mathsf{m}%
\right\vert ^{s}\right) ^{k}\left\vert p_{\mathsf{m}}\right\vert ^{s}<+\infty
\label{convergencepremiere}
\end{equation}%
for $s\in \left[ 1,+\infty \right) $, $w=$ $\left( w_{\mathsf{m}}>0\right) _{%
\mathsf{m}\in \mathbb{Z}}$ a sequence of positive weights and $k\in \mathbb{R%
}$, for which various compact embedding theorems were proven. Relation (\ref%
{convergencepremiere}) was motivated in the first place by the usual Sobolev
space theory of periodic functions and its relation to Fourier analysis,
whereby $s$ was referred to as the degree of summability of $\mathsf{p}$ and 
$k$ to its generalized order of differentiability. As a matter of fact, the $%
h_{\mathbb{C},w}^{k,s}$ in turn were introduced as a generalization of the
Hilbert space $h_{\mathbb{C},w}^{1,2}$ that played a decisive role in the
analysis of a class of master equations with non-constant coefficients
arising in non-equilibrium statistical mechanics, which was carried out in 
\cite{vuillermotbis}. Accordingly, the remaining part of this article is
organized as follows: In Section 2 we generalize (\ref{convergencepremiere})
in a natural way by defining an Orlicz-Sobolev class $\mathsf{c}_{k,\Phi ,w}$
from a continuous and convex function $\Phi $ that is not limited to
polynomial growth and which satisfies additional properties. The price to
pay for such a generalization is that $\mathsf{c}_{k,\Phi ,w}$ is no longer
a vector space, though it remains a convex set. As a consequence we then
define two \textit{bona fide} Banach sequence spaces around $\mathsf{c}%
_{k,\Phi ,w}$, namely, the smallest subspace containing $\mathsf{c}_{k,\Phi
,w}$ on the one hand, and the largest subspace contained in $\mathsf{c}%
_{k,\Phi ,w}$ on the other hand, denoted by $l_{k,\Phi ,w}$ and $\mathsf{s}%
_{k,\Phi ,w}$ and referred to as the large Orlicz-Sobolev sequence space and
the small Orlicz-Sobolev sequence space, respectively. We do so in a way
that allows us to establish the equivalence with other definitions of such
spaces when $k=0$ and $w_{\mathsf{m}}=1$ for every $\mathsf{m}$, for example
those set forth at the beginning of Chapter 4 in \cite{lindenstrausstra}. We
proceed by proving the existence of continuous embeddings for $l_{k,\Phi ,w}$
and $\mathsf{s}_{k,\Phi ,w}$ and by constructing a Schauder basis in $%
\mathsf{s}_{k,\Phi ,w}$. Using that basis in an essential way we then prove
the existence of various compact embeddings for $\mathsf{s}_{k,\Phi ,w}$
under the condition that $\mathsf{s}_{k,\Phi ,w}=l_{k,\Phi ,w}$, an equality
guaranteed by imposing a constraint on $\Phi $ regarding its behavior in a
small neighborhood of the origin. We devote Section 3 to illustrating the
results of Section 2 with concrete examples.

As far as we know the embedding theorems stated below are new and provide
complementary information regarding the results stated in Section 10.3 of
Chapter 10 in \cite{raoren2} we alluded to above. Moreover, as amply
testified by the many references appearing in \cite{raoren1} and \cite%
{raoren2} to which we refer the reader for more detailed information, Orlicz
sequence spaces play a crucial role in a number of important applications
including very recent analyses of discrete Wiener-Hopf operators as in \cite%
{karlovychthampi}. It is therefore also our contention that the new scale we
introduce will be as useful for applications to various areas of mathematics
or mathematical physics as the spaces $h_{\mathbb{C},w}^{k,s}$ are (see,
e.g., \cite{vuillermotbis}).

\section{Main results}

Let $\Phi :\left[ 0,+\infty \right) \rightarrow \left[ 0,+\infty \right) $
be a continuous non-decreasing and convex function satisfying $\Phi \left(
0\right) =0$, $\Phi \left( t\right) >0$ for every $t>0$ and $%
\lim_{t\rightarrow +\infty }\Phi \left( t\right) =+\infty $, that is, a 
\textit{non-degenerate Orlicz function} in the sense of Definition 4.a.1. in
Chapter 4 of \cite{lindenstrausstra}. Under these conditions, it is plain
that $\Phi $ is actually strictly increasing throughout the interval $\left[
0,+\infty \right) $. Henceforth any function that satisfies the above
properties will be simply referred to a an \textit{Orlicz function}. Let $w=$
$\left( w_{\mathsf{m}}\right) _{\mathsf{m}\in \mathbb{Z}}$ be a sequence of
weights satisfying $w_{\mathsf{m}}>0$ for every $\mathsf{m}$ and let $k\in 
\mathbb{R}$. We begin with the following:

\bigskip

\textbf{Definition.} We write $\mathsf{c}_{k,\Phi ,w}$ for the \textit{%
weighted Orlicz-Sobolev sequence class} consisting of all complex sequences $%
\mathsf{p}=\left( p_{\mathsf{m}}\right) _{\mathsf{m}\in \mathbb{Z}}$
satisfying%
\begin{equation}
\sum_{\mathsf{m}\in \mathbb{Z}}w_{\mathsf{m}}\left( 1+\Phi \left( \left\vert 
\mathsf{m}\right\vert \right) \right) ^{k}\Phi \left( \left\vert p_{\mathsf{m%
}}\right\vert \right) <+\infty .  \label{convergence}
\end{equation}%
For the particular case $\Phi \left( t\right) =t^{s}$ with $s\in \left[
1,+\infty \right) $, the preceding expression leads to%
\begin{equation}
\left\Vert \mathsf{p}\right\Vert _{k,s,w}:=\left( \sum_{\mathsf{m}\in 
\mathbb{Z}}w_{\mathsf{m}}\left( 1+\left\vert \mathsf{m}\right\vert
^{s}\right) ^{k}\left\vert p_{\mathsf{m}}\right\vert ^{s}\right) ^{\frac{1}{s%
}}<+\infty ,  \label{norm}
\end{equation}%
which is the norm carried by the Banach spaces $h_{\mathbb{C},w}^{k,s}$ we
alluded to in the introduction. In case of a general Orlicz function,
however, expression (\ref{convergence}) does not define a norm and the class 
$\mathsf{c}_{k,\Phi ,w}$ does not inherit a linear structure unless further
conditions are imposed on $\Phi $, as is already well known for the case $%
k=0 $ and $w_{\mathsf{m}}=1$ for every $\mathsf{m}$. It remains nonetheless
a convex set as a consequence of the monotonicity and the convexity of $\Phi 
$. In the present context this means that if $\mathsf{p,q}\in \mathsf{c}%
_{k,\Phi ,w}$, then $\lambda \mathsf{p}+\mu \mathsf{q}\in \mathsf{c}_{k,\Phi
,w}$ for all $\lambda ,\mu \in \mathbb{C}$ satisfying $\left\vert \lambda
\right\vert +\left\vert \mu \right\vert \leq 1$. Therefore, as is the case
for Orlicz classes defined with respect to Lebesgue measure in Euclidean
space as in \cite{adamsfournier}, or with respect to more general diffusive
measures as in \cite{raoren1} or \cite{raoren2}, we now construct two 
\textit{bona fide }Banach spaces from $\mathsf{c}_{k,\Phi ,w}$ and establish
relations between them. In order to do so we write $\mu _{k,\Phi ,w}=\left(
\mu _{k,\Phi ,w,\mathsf{m}}\right) _{\mathsf{m}\in \mathbb{Z}}$ where 
\begin{equation}
\mu _{k,\Phi ,w,\mathsf{m}}:=w_{\mathsf{m}}\left( 1+\Phi \left( \left\vert 
\mathsf{m}\right\vert \right) \right) ^{k},  \label{simplification}
\end{equation}%
and only use the explicit form of the right-hand side wherever it is
relevant. Our starting point is the following elementary result:

\bigskip

\textbf{Proposition 1.} \textit{For every} $\mathsf{p}\in \mathsf{c}_{k,\Phi
,w}$ \textit{we have}%
\begin{equation}
\left\Vert \mathsf{p}\right\Vert _{k,\Phi ,w}:=\inf \left\{ \rho >0:\sum_{%
\mathsf{m}\in \mathbb{Z}}\mu _{k,\Phi ,w,\mathsf{m}}\Phi \left( \frac{%
\left\vert p_{\mathsf{m}}\right\vert }{\rho }\right) \leq 1\right\} <+\infty
.  \label{luxnorm}
\end{equation}

\bigskip

\textbf{Proof. }The function%
\begin{equation*}
\rho \rightarrow \sum_{\mathsf{m}\in \mathbb{Z}}\mu _{k,\Phi ,w,\mathsf{m}%
}\Phi \left( \frac{\left\vert p_{\mathsf{m}}\right\vert }{\rho }\right)
\end{equation*}%
remains finite and non-increasing for $\rho \in \left[ 1,+\infty \right) $.
Moreover, if $\left( \rho _{\mathsf{n}}\right) _{\mathsf{n}\in \mathbb{N}}$
is any sequence satisfying $\rho _{\mathsf{n}}>1$ with $\rho _{\mathsf{n}%
}\rightarrow +\infty $ as $\mathsf{n}\rightarrow +\infty $ we have%
\begin{equation*}
\mu _{k,\Phi ,w,\mathsf{m}}\Phi \left( \frac{\left\vert p_{\mathsf{m}%
}\right\vert }{\rho _{\mathsf{n}}}\right) \leq \mu _{k,\Phi ,w,\mathsf{m}%
}\Phi \left( \left\vert p_{\mathsf{m}}\right\vert \right)
\end{equation*}%
and%
\begin{equation*}
\lim_{\mathsf{n}\rightarrow +\infty }\mu _{k,\Phi ,w,\mathsf{m}}\Phi \left( 
\frac{\left\vert p_{\mathsf{m}}\right\vert }{\rho _{\mathsf{n}}}\right) =0
\end{equation*}%
for every $\mathsf{m}\in \mathbb{Z}$, so that%
\begin{equation*}
\lim_{\mathsf{n}\rightarrow +\infty }\sum_{\mathsf{m}\in \mathbb{Z}}\mu
_{k,\Phi ,w,\mathsf{m}}\Phi \left( \frac{\left\vert p_{\mathsf{m}%
}\right\vert }{\rho _{\mathsf{n}}}\right) =0
\end{equation*}%
as well by dominated convergence since (\ref{convergence}) holds. Therefore,
for all sufficiently large $\rho $ we have in particular%
\begin{equation*}
\sum_{\mathsf{m}\in \mathbb{Z}}\mu _{k,\Phi ,w,\mathsf{m}}\Phi \left( \frac{%
\left\vert p_{\mathsf{m}}\right\vert }{\rho }\right) \leq 1,
\end{equation*}%
thereby proving the result. \ \ $\blacksquare $

\bigskip

\textsc{Remarks.} (1) If $\Phi \left( t\right) =t^{s}$ with $s\in \left[
1,+\infty \right) $, it is plain that (\ref{luxnorm}) reduces to (\ref{norm}%
).

(2) We have%
\begin{equation}
\sum_{\mathsf{m}\in \mathbb{Z}}\mu _{k,\Phi ,w,\mathsf{m}}\Phi \left( \frac{%
\left\vert p_{\mathsf{m}}\right\vert }{\left\Vert \mathsf{p}\right\Vert
_{k,\Phi ,w}}\right) \leq 1  \label{crucialrelation}
\end{equation}%
whenever $\left\Vert \mathsf{p}\right\Vert _{k,\Phi ,w}>0$, as a consequence
of the monotone convergence theorem. To see that it is sufficient to choose
a minimizing sequence $\left( \rho _{\mathsf{n}}\right) _{\mathsf{n}\in 
\mathbb{N}}$ satisfying $\rho _{\mathsf{n+1}}<\rho _{\mathsf{n}}$ and 
\begin{equation*}
\sum_{\mathsf{m}\in \mathbb{Z}}\mu _{k,\Phi ,w,\mathsf{m}}\Phi \left( \frac{%
\left\vert p_{\mathsf{m}}\right\vert }{\rho _{\mathsf{n}}}\right) \leq 1
\end{equation*}%
for every $\mathsf{n}\in \mathbb{N}$, with $\rho _{\mathsf{n}}\rightarrow
\left\Vert \mathsf{p}\right\Vert _{k,\Phi ,w}$ as $\mathsf{n\rightarrow
+\infty }$. Relation (\ref{crucialrelation}) plays a crucial role in what
follows.

\bigskip

Now let $l_{k,\Phi ,w}$ be the set of \textit{all} complex sequences such
that (\ref{luxnorm}) holds. It is easily seen that $l_{k,\Phi ,w}$ is a
normed vector space with respect to the usual algebraic operations and (\ref%
{luxnorm}), which indeed possesses all the attributes of a norm. In
particular, the triangle inequality is a direct consequence of the
monotonicity and the convexity of $\Phi $ together with (\ref%
{crucialrelation}). For if $\left\Vert \mathsf{p}\right\Vert _{k,\Phi ,w}>0$
and $\left\Vert \mathsf{q}\right\Vert _{k,\Phi ,w}>0$, we have successively%
\begin{eqnarray*}
&&\sum_{\mathsf{m}\in \mathbb{Z}}\mu _{k,\Phi ,w,\mathsf{m}}\Phi \left( 
\frac{\left\vert p_{\mathsf{m}}+q_{\mathsf{m}}\right\vert }{\left\Vert 
\mathsf{p}\right\Vert _{k,\Phi ,w}+\left\Vert \mathsf{q}\right\Vert _{k,\Phi
,w}}\right) \\
&\leq &\frac{\left\Vert \mathsf{p}\right\Vert _{k,\Phi ,w}}{\left\Vert 
\mathsf{p}\right\Vert _{k,\Phi ,w}+\left\Vert \mathsf{q}\right\Vert _{k,\Phi
,w}}\sum_{\mathsf{m}\in \mathbb{Z}}\mu _{k,\Phi ,w,\mathsf{m}}\Phi \left( 
\frac{\left\vert p_{\mathsf{m}}\right\vert }{\left\Vert \mathsf{p}%
\right\Vert _{k,\Phi ,w}}\right) \\
&&+\frac{\left\Vert \mathsf{q}\right\Vert _{k,\Phi ,w}}{\left\Vert \mathsf{p}%
\right\Vert _{k,\Phi ,w}+\left\Vert \mathsf{q}\right\Vert _{k,\Phi ,w}}\sum_{%
\mathsf{m}\in \mathbb{Z}}\mu _{k,\Phi ,w,\mathsf{m}}\Phi \left( \frac{%
\left\vert q_{\mathsf{m}}\right\vert }{\left\Vert \mathsf{q}\right\Vert
_{k,\Phi ,w}}\right) \\
&\leq &\frac{\left\Vert \mathsf{p}\right\Vert _{k,\Phi ,w}}{\left\Vert 
\mathsf{p}\right\Vert _{k,\Phi ,w}+\left\Vert \mathsf{q}\right\Vert _{k,\Phi
,w}}+\frac{\left\Vert \mathsf{q}\right\Vert _{k,\Phi ,w}}{\left\Vert \mathsf{%
p}\right\Vert _{k,\Phi ,w}+\left\Vert \mathsf{q}\right\Vert _{k,\Phi ,w}}=1,
\end{eqnarray*}%
which immediately implies that $\left\Vert \mathsf{p+q}\right\Vert _{k,\Phi
,w}\leq \left\Vert \mathsf{p}\right\Vert _{k,\Phi ,w}+\left\Vert \mathsf{q}%
\right\Vert _{k,\Phi ,w}<+\infty $, this inequality being obvious if $%
\left\Vert \mathsf{p}\right\Vert _{k,\Phi ,w}=0$ or if $\left\Vert \mathsf{q}%
\right\Vert _{k,\Phi ,w}=0$. Consequently, from the preceding considerations
we infer that%
\begin{equation}
\limfunc{span}\mathsf{c}_{k,\Phi ,w}\subseteq l_{k,\Phi ,w}.
\label{inclusion1}
\end{equation}%
In fact, the following result holds:

\bigskip

\textbf{Theorem 1.} \textit{We have}%
\begin{equation}
\limfunc{span}\mathsf{c}_{k,\Phi ,w}=l_{k,\Phi ,w}  \label{equality}
\end{equation}%
\textit{and }$l_{k,\Phi ,w}$ \textit{is a complex Banach space with respect
to (\ref{luxnorm}).}

\bigskip

\textbf{Proof.} For every $\mathsf{p}\in $ $l_{k,\Phi ,w}$ there exists $%
\rho >0$ such that 
\begin{equation*}
\sum_{\mathsf{m}\in \mathbb{Z}}\mu _{k,\Phi ,w,\mathsf{m}}\Phi \left( \frac{%
\left\vert p_{\mathsf{m}}\right\vert }{\rho }\right) \leq 1
\end{equation*}%
according to (\ref{luxnorm}), so that by choosing $\mathsf{q}=\rho ^{-1}%
\mathsf{p}$ we have $\mathsf{q}\in $ $\mathsf{c}_{k,\Phi ,w}$ with $\mathsf{%
p=}\rho \mathsf{q}$. Consequently%
\begin{equation*}
l_{k,\Phi ,w}\subseteq \limfunc{span}\mathsf{c}_{k,\Phi ,w},
\end{equation*}%
which proves (\ref{equality}) since (\ref{inclusion1}) holds. The proof that 
$l_{k,\Phi ,w}$ is complete with respect to (\ref{luxnorm}) results from a
modification of standard arguments and is therefore omitted (see, e.g., \cite%
{adamsfournier}, \cite{krasnorut}, \cite{raoren1} and \cite{raoren2}). \ \ $%
\blacksquare $

\bigskip

\textsc{Remark.} Our definition of $l_{k,\Phi ,w}$ implies that for every $%
\mathsf{p}\in l_{k,\Phi ,w}$, there exists at least one $\rho >0$ such that%
\begin{equation}
\sum_{\mathsf{m}\in \mathbb{Z}}\mu _{k,\Phi ,w,\mathsf{m}}\Phi \left( \frac{%
\left\vert p_{\mathsf{m}}\right\vert }{\rho }\right) <+\infty .
\label{convergencebis}
\end{equation}%
\newline
Conversely, if this property holds true then, by using a dominated
convergence argument similar to that set forth in the proof of Proposition
1, for all sufficiently large $\rho >0$ we get%
\begin{equation*}
\sum_{\mathsf{m}\in \mathbb{Z}}\mu _{k,\Phi ,w,\mathsf{m}}\Phi \left( \frac{%
\left\vert p_{\mathsf{m}}\right\vert }{\rho }\right) \leq 1
\end{equation*}%
and hence $\left\Vert \mathsf{p}\right\Vert _{k,\Phi ,w}<+\infty .$
Therefore, we might have defined $l_{k,\Phi ,w}$ as the space consisting of
all complex sequences $\mathsf{p}$ for which there exists \textit{at least
one }$\rho >0$ depending on $\mathsf{p}$ implying (\ref{convergencebis}).
With $k=0$ and $w_{\mathsf{m}}=1$ for every $\mathsf{m}\in \mathbb{Z}$, this
establishes the connection with the definition of the Orlicz sequence space
given at the beginning of Chapter 4 in \cite{lindenstrausstra}.

\bigskip

In what follows we refer to $l_{k,\Phi ,w}$ as the\textit{\ large
Orlicz-Sobolev sequence space,} as we can also define the \textit{small
Orlicz-Sobolev sequence} \textit{space} $\mathsf{s}_{k,\Phi ,w}$ consisting
of all complex sequences for which (\ref{convergencebis}) holds for \textit{%
every }$\rho >0$. Again by virtue of the monotonicity and the convexity of $%
\Phi $, it follows that $\mathsf{s}_{k,\Phi ,w}$ is a vector subspace of $%
l_{k,\Phi ,w}$. Actually the following result holds, which justifies \textit{%
a posteriori} the chosen terminology:

\bigskip

\textbf{Theorem 2.} \textit{The linear space }$\mathsf{s}_{k,\Phi ,w}$%
\textit{\ is a complex Banach space, and is the largest subspace contained
in the weighted Orlicz-Sobolev sequence class }$\mathsf{c}_{k,\Phi ,w}$. 
\textit{Moreover, }$\mathsf{s}_{k,\Phi ,w}$\textit{\ is separable.}

\textit{\bigskip }

\textbf{Proof.} Let $\left( \mathsf{p}_{\mathsf{N}}\right) $ be a sequence
in $\mathsf{s}_{k,\Phi ,w}$ converging to $\mathsf{p}$ with respect to (\ref%
{luxnorm}). In order to prove the first part of the proposition, it is
sufficient to show that $\mathsf{p}\in \mathsf{s}_{k,\Phi ,w}$. Writing $%
\mathsf{p=p}_{\mathsf{N}}+\left( \mathsf{p-p}_{\mathsf{N}}\right) $ and
owing once again to the monotonicity and the convexity of $\Phi $, we have
successively%
\begin{eqnarray*}
&&\sum_{\mathsf{m}\in \mathbb{Z}}\mu _{k,\Phi ,w,\mathsf{m}}\Phi \left( 
\frac{\left\vert p_{\mathsf{m}}\right\vert }{2\rho }\right) \\
&\leq &\sum_{\mathsf{m}\in \mathbb{Z}}\mu _{k,\Phi ,w,\mathsf{m}}\Phi \left( 
\frac{\left\vert \left( \mathsf{p}_{\mathsf{N}}\right) _{\mathsf{m}%
}\right\vert }{2\rho }+\frac{\left\vert p_{\mathsf{m}}-\left( \mathsf{p}_{%
\mathsf{N}}\right) _{\mathsf{m}}\right\vert }{2\rho }\right) \\
&\leq &\frac{1}{2}\sum_{\mathsf{m}\in \mathbb{Z}}\mu _{k,\Phi ,w,\mathsf{m}%
}\Phi \left( \frac{\left\vert \left( \mathsf{p}_{\mathsf{N}}\right) _{%
\mathsf{m}}\right\vert }{\rho }\right) +\frac{1}{2}\sum_{\mathsf{m}\in 
\mathbb{Z}}\mu _{k,\Phi ,w,\mathsf{m}}\Phi \left( \frac{\left\vert p_{%
\mathsf{m}}-\left( \mathsf{p}_{\mathsf{N}}\right) _{\mathsf{m}}\right\vert }{%
\rho }\right) .
\end{eqnarray*}%
The first series on the right hand-side is convergent for every $\rho >0$
since $\mathsf{p}_{\mathsf{N}}\in \mathsf{s}_{k,\Phi ,w}$. We now show that
the same property holds true for the second series by virtue of the fact
that $\left\Vert \mathsf{p-p}_{\mathsf{N}}\right\Vert _{k,\Phi
,w}\rightarrow 0$ as $\mathsf{N\rightarrow +\infty }$. Indeed, for each $%
\epsilon >0$ there exists $\mathsf{N(}\epsilon )\in $ $\mathbb{N}^{+}$ such
that 
\begin{equation}
\inf \left\{ \rho >0:\sum_{\mathsf{m}\in \mathbb{Z}}\mu _{k,\Phi ,w,\mathsf{m%
}}\Phi \left( \frac{\left\vert \left( \mathsf{p-p}_{\mathsf{N}}\right) _{%
\mathsf{m}}\right\vert }{\rho }\right) \leq 1\right\} \leq \epsilon
\label{inequality2}
\end{equation}%
for every $\mathsf{N\geq N(}\epsilon )$. For each such $\mathsf{N}$ it is
therefore sufficient to prove that%
\begin{equation}
\sum_{\mathsf{m}\in \mathbb{Z}}\mu _{k,\Phi ,w,\mathsf{m}}\Phi \left( \frac{%
\left\vert \left( \mathsf{p-p}_{\mathsf{N}}\right) _{\mathsf{m}}\right\vert 
}{\rho }\right) \leq 1  \label{inequality3}
\end{equation}%
for every $\rho >0$. Let us suppose that the preceding inequality does not
hold. Then there exists $\rho ^{\ast }>0$ such that%
\begin{equation*}
\sum_{\mathsf{m}\in \mathbb{Z}}\mu _{k,\Phi ,w,\mathsf{m}}\Phi \left( \frac{%
\left\vert \left( \mathsf{p-p}_{\mathsf{N}}\right) _{\mathsf{m}}\right\vert 
}{\rho ^{\ast }}\right) >1,
\end{equation*}%
and the same estimate holds for every $\rho \in \left( 0,\rho ^{\ast
}\right) $ by virtue of the properties of $\Phi $. Consequently, we
necessarily have the inclusion%
\begin{equation*}
\left\{ \rho >0:\sum_{\mathsf{m}\in \mathbb{Z}}\mu _{k,\Phi ,w,\mathsf{m}%
}\Phi \left( \frac{\left\vert \left( \mathsf{p-p}_{\mathsf{N}}\right) _{%
\mathsf{m}}\right\vert }{\rho }\right) \leq 1\right\} \subseteq \left( \rho
^{\ast },+\infty \right)
\end{equation*}%
and therefore%
\begin{equation*}
\inf \left( \rho ^{\ast },+\infty \right) \leq \epsilon
\end{equation*}%
for every $\epsilon >0$ according to (\ref{inequality2}), which is
impossible if $\epsilon \in \left( 0,\rho ^{\ast }\right) $. Consequently (%
\ref{inequality3}) holds, and hence\textit{\ a fortiori}%
\begin{equation*}
\sum_{\mathsf{m}\in \mathbb{Z}}\mu _{k,\Phi ,w,\mathsf{m}}\Phi \left( \frac{%
\left\vert p_{\mathsf{m}}\right\vert }{2\rho }\right) <+\infty
\end{equation*}%
for every $\rho >0$ so that indeed $\mathsf{p}\in \mathsf{s}_{k,\Phi ,w}$.

While it is plain that $\mathsf{s}_{k,\Phi ,w}\subseteq \mathsf{c}_{k,\Phi
,w}$, we now show that if $\mathsf{S}$ is a linear subspace of $\mathsf{c}%
_{k,\Phi ,w}$ then $\mathsf{S}\subseteq \mathsf{s}_{k,\Phi ,w}$. Let $\rho
>0 $ be arbitrary and let us choose $\lambda \in \mathbb{C}$ such that $%
\left\vert \lambda \right\vert ^{-1}=\rho $. For any $\mathsf{p\in S}$ we
then have $\lambda \mathsf{p}\in \mathsf{S}\subseteq $ $\mathsf{c}_{k,\Phi
,w}$ and thereby%
\begin{eqnarray*}
&&\sum_{\mathsf{m}\in \mathbb{Z}}\mu _{k,\Phi ,w,\mathsf{m}}\Phi \left( 
\frac{\left\vert p_{\mathsf{m}}\right\vert }{\rho }\right) \\
&=&\sum_{\mathsf{m}\in \mathbb{Z}}\mu _{k,\Phi ,w,\mathsf{m}}\Phi \left( 
\frac{\left\vert p_{\mathsf{m}}\right\vert }{\left\vert \lambda \right\vert
^{-1}}\right) =\sum_{\mathsf{m}\in \mathbb{Z}}\mu _{k,\Phi ,w,\mathsf{m}%
}\Phi \left( \left\vert \left( \lambda \mathsf{p}\right) _{\mathsf{m}%
}\right\vert \right) <+\infty ,
\end{eqnarray*}%
which is the desired result.

As for separability, we proceed by showing that the set of sequences $%
\mathsf{e}_{\mathsf{m}}$ defined by $\left( \mathsf{e}_{\mathsf{m}}\right) _{%
\mathsf{n}}=\delta _{\mathsf{m,n}}$ for all $\mathsf{m,n}\in \mathbb{Z}$
constitutes a Schauder basis of $\mathsf{s}_{k,\Phi ,w}$ (we refer the
reader for instance to \cite{albiackalton} and \cite{lindenstrausstra} for
basic definitions and many properties regarding such bases). This means in
particular that every $\mathsf{p=}\left( p_{\mathsf{m}}\right) \in \mathsf{s}%
_{k,\Phi ,w}$ may be expanded in a unique fashion as the norm-convergent
series 
\begin{equation*}
\mathsf{p=}\sum_{\mathsf{m}\in \mathbb{Z}}p_{\mathsf{m}}\mathsf{e}_{\mathsf{m%
}}\text{.}
\end{equation*}%
First, it is plain that $\mathsf{e}_{\mathsf{m}}\in \mathsf{s}_{k,\Phi ,w}$
for each $\mathsf{m}$, so that for every $\mathsf{M}\in \mathbb{N}^{+}$ we
may define%
\begin{equation}
\mathsf{p}_{\mathsf{M}}\mathsf{=}\sum_{\left\vert \mathsf{m}\right\vert \leq 
\mathsf{M}}p_{\mathsf{m}}\mathsf{e}_{\mathsf{m}}  \label{cauchysequence}
\end{equation}%
as an element of that space. It is then sufficient to show that (\ref%
{cauchysequence}) generates a Cauchy sequence therein. That is, for each $%
\epsilon >0$ we have to prove the existence of an $\mathsf{N}\left( \epsilon
\right) \in \mathbb{N}^{+}$ such that%
\begin{equation}
\sum_{\mathsf{m}\in \mathbb{Z}}\mu _{k,\Phi ,w,\mathsf{m}}\Phi \left( \frac{%
\left\vert \left( \mathsf{p}_{\mathsf{M}}-\mathsf{p}_{\mathsf{N}}\right) _{%
\mathsf{m}}\right\vert }{\epsilon }\right) \leq 1  \label{cauchy}
\end{equation}%
whenever $\mathsf{M,N\geq N}\left( \epsilon \right) $, for then the relation 
\begin{equation*}
\left\Vert \mathsf{p}_{\mathsf{M}}-\mathsf{p}_{\mathsf{N}}\right\Vert
_{k,\Phi ,w}\leq \epsilon
\end{equation*}%
follows from (\ref{luxnorm}). For any $\mathsf{N}\in \mathbb{N}^{+}$ with $%
\mathsf{M>N}$ we have%
\begin{equation*}
\mathsf{p}_{\mathsf{M}}-\mathsf{p}_{\mathsf{N}}=\sum_{\mathsf{m}=-\mathsf{M}%
}^{-\left( \mathsf{N}+1\right) }p_{\mathsf{m}}\mathsf{e}_{\mathsf{m}}+\sum_{%
\mathsf{m}=\mathsf{N+}1}^{\mathsf{M}}p_{\mathsf{m}}\mathsf{e}_{\mathsf{m}}
\end{equation*}%
and thereby%
\begin{equation*}
\left( \mathsf{p}_{\mathsf{M}}-\mathsf{p}_{\mathsf{N}}\right) _{\mathsf{m}%
}=\left\{ 
\begin{array}{c}
p_{\mathsf{m}}\text{ if }\mathsf{m}\in \mathbb{Z}\cap \left[ -\mathsf{M,}%
-\left( \mathsf{N}+1\right) \right] \text{ or if }\mathsf{m}\in \text{ }%
\mathbb{Z}\cap \left[ \mathsf{N}+1\mathsf{,M}\right] , \\ 
\\ 
0\text{ otherwise}%
\end{array}%
\right.
\end{equation*}%
for the $\mathsf{m}^{\mathsf{th}}$ component of $\mathsf{p}_{\mathsf{M}}-%
\mathsf{p}_{\mathsf{N}}$. Therefore we get%
\begin{eqnarray}
&&\sum_{\mathsf{m}\in \mathbb{Z}}\mu _{k,\Phi ,w,\mathsf{m}}\Phi \left( 
\frac{\left\vert \left( \mathsf{p}_{\mathsf{M}}-\mathsf{p}_{\mathsf{N}%
}\right) _{\mathsf{m}}\right\vert }{\epsilon }\right)  \notag \\
&=&\sum_{\mathsf{m}=-\mathsf{M}}^{-\left( \mathsf{N}+1\right) }\mu _{k,\Phi
,w,\mathsf{m}}\Phi \left( \frac{\left\vert p_{\mathsf{m}}\right\vert }{%
\epsilon }\right) +\sum_{\mathsf{m}=\mathsf{N}+1}^{\mathsf{M}}\mu _{k,\Phi
,w,\mathsf{m}}\Phi \left( \frac{\left\vert p_{\mathsf{m}}\right\vert }{%
\epsilon }\right) .  \label{equalitybis}
\end{eqnarray}%
Furthermore, since $\mathsf{p}\in \mathsf{s}_{k,\Phi ,w}$ we have%
\begin{equation}
\sum_{\mathsf{m}\in \mathbb{Z}}\mu _{k,\Phi ,w,\mathsf{m}}\Phi \left( \frac{%
\left\vert p_{\mathsf{m}}\right\vert }{\epsilon }\right) =\lim_{\mathsf{%
M\rightarrow +\infty }}\sum_{\left\vert \mathsf{m}\right\vert \leq \mathsf{M}%
}\mu _{k,\Phi ,w,\mathsf{m}}\Phi \left( \frac{\left\vert p_{\mathsf{m}%
}\right\vert }{\epsilon }\right) <+\infty ,  \label{convergencequinto}
\end{equation}%
which implies that%
\begin{equation*}
\lim_{\mathsf{M\rightarrow +\infty }}\sum_{\mathsf{m=-M}}^{\mathsf{0}}\mu
_{k,\Phi ,w,\mathsf{m}}\Phi \left( \frac{\left\vert p_{\mathsf{m}%
}\right\vert }{\epsilon }\right) <+\infty
\end{equation*}%
and%
\begin{equation*}
\lim_{\mathsf{M\rightarrow +\infty }}\sum_{\mathsf{m=0}}^{\mathsf{M}}\mu
_{k,\Phi ,w,\mathsf{m}}\Phi \left( \frac{\left\vert p_{\mathsf{m}%
}\right\vert }{\epsilon }\right) <+\infty
\end{equation*}%
since all terms in the above sums are non-negative. Thus, the mappings%
\begin{equation*}
\mathsf{M\rightarrow }\sum_{\mathsf{m=-M}}^{\mathsf{0}}\mu _{k,\Phi ,w,%
\mathsf{m}}\Phi \left( \frac{\left\vert p_{\mathsf{m}}\right\vert }{\epsilon 
}\right)
\end{equation*}%
and 
\begin{equation*}
\mathsf{M\rightarrow }\sum_{\mathsf{m=0}}^{\mathsf{M}}\mu _{k,\Phi ,w,%
\mathsf{m}}\Phi \left( \frac{\left\vert p_{\mathsf{m}}\right\vert }{\epsilon 
}\right)
\end{equation*}%
are Cauchy sequences in $\mathbb{R}$. Therefore, there exists $\mathsf{N}%
_{1}\left( \epsilon \right) \in \mathbb{N}^{+}$ such that%
\begin{equation*}
\sum_{\mathsf{m=-M}}^{\mathsf{-}\left( \mathsf{N}+1\right) }\mu _{k,\Phi ,w,%
\mathsf{m}}\Phi \left( \frac{\left\vert p_{\mathsf{m}}\right\vert }{\epsilon 
}\right) \leq \frac{1}{2}
\end{equation*}%
for all $\mathsf{M>N\geq N}_{1}\left( \epsilon \right) $, and similarly
there exists $\mathsf{N}_{2}\left( \epsilon \right) \in \mathbb{N}^{+}$ such
that%
\begin{equation*}
\sum_{\mathsf{m=N+1}}^{\mathsf{M}}\mu _{k,\Phi ,w,\mathsf{m}}\Phi \left( 
\frac{\left\vert p_{\mathsf{m}}\right\vert }{\epsilon }\right) \leq \frac{1}{%
2}
\end{equation*}%
whenever $\mathsf{M>N\geq N}_{2}\left( \epsilon \right) $. Consequently,
choosing $\mathsf{N}\left( \epsilon \right) =\mathsf{N}_{1}\left( \epsilon
\right) \vee \mathsf{N}_{2}\left( \epsilon \right) $ and taking $\mathsf{%
M>N\geq N}\left( \epsilon \right) $ we obtain (\ref{cauchy}) according to (%
\ref{equalitybis}), and thereby%
\begin{equation*}
\mathsf{p=}\lim_{\mathsf{M\rightarrow +\infty }}\sum_{\left\vert \mathsf{m}%
\right\vert \leq \mathsf{M}}p_{\mathsf{m}}\mathsf{e}_{\mathsf{m}}=\sum_{%
\mathsf{m}\in \mathbb{Z}}p_{\mathsf{m}}\mathsf{e}_{\mathsf{m}}
\end{equation*}%
which is the desired result. \ \ $\blacksquare $

\bigskip

\textsc{Remark.} It is essential that $\mathsf{p}\in \mathsf{s}_{k,\Phi ,w}$
for (\ref{convergencequinto}) to hold for \textit{every} $\epsilon >0$. If $%
\mathsf{p}\in l_{k,\Phi ,w}$ but $\mathsf{p}\notin \mathsf{s}_{k,\Phi ,w}$
then the above argument breaks down, and indeed the sequence $\left( \mathsf{%
e}_{\mathsf{m}}\right) \mathsf{\ }$\textit{does not} constitute a Schauder
basis of $l_{k,\Phi ,w}$ nor is this space separable in general, as already
pointed out and proven in a different way in \cite{lindenstrausstra} for
what corresponds to the case $k=0$ and $w_{\mathsf{m}}=1$ for every $\mathsf{%
m}$ in this article (see, e.g., Proposition 4.a.4 in that reference).

\bigskip

We have already pointed to the fact that the Orlicz-Sobolev sequence class $%
\mathsf{c}_{k,\Phi ,w}$ does not carry a linear structure in general, so
that in such a case Theorems 1 and 2 imply that $\mathsf{c}_{k,\Phi ,w}$ and
the spaces $\mathsf{s}_{k,\Phi ,w}$ and $l_{k,\Phi ,w}$ are all distinct.
However, at this stage it is already possible to prove continuous embedding
results for $l_{k,\Phi ,w}$ and $\mathsf{s}_{k,\Phi ,w}$ in a relatively
simple way, in sharp contrast to the compact case as we shall see below. As
a preliminary remark we note that the embedding

\begin{equation*}
l_{k^{\prime },\Phi ,w}\rightarrow l_{k,\Phi ,w}
\end{equation*}%
is continuous when $k^{\prime }\geq k$ by virtue of (\ref{luxnorm}) and the
obvious inclusion%
\begin{eqnarray*}
&&\left\{ \rho >0:\sum_{\mathsf{m}\in \mathbb{Z}}w_{\mathsf{m}}\left( 1+\Phi
\left( \left\vert \mathsf{m}\right\vert \right) \right) ^{k^{\prime }}\Phi
\left( \frac{\left\vert p_{\mathsf{m}}\right\vert }{\rho }\right) \leq
1\right\} \\
&\subseteq &\left\{ \rho >0:\sum_{\mathsf{m}\in \mathbb{Z}}w_{\mathsf{m}%
}\left( 1+\Phi \left( \left\vert \mathsf{m}\right\vert \right) \right)
^{k}\Phi \left( \frac{\left\vert p_{\mathsf{m}}\right\vert }{\rho }\right)
\leq 1\right\} .
\end{eqnarray*}%
Similarly we have 
\begin{equation}
\mathsf{s}_{k^{\prime },\Phi ,w}\rightarrow \mathsf{s}_{k,\Phi ,w}.
\label{embeddingbis}
\end{equation}%
Indeed, in order to get (\ref{embeddingbis}) it is sufficient to show that $%
\mathsf{s}_{k^{\prime },\Phi ,w}\subseteq \mathsf{s}_{k,\Phi ,w}$ since $%
\mathsf{s}_{k^{\prime },\Phi ,w}\subseteq l_{k^{\prime },\Phi ,w}$. But for
each $\mathsf{p}\in \mathsf{s}_{k^{\prime },\Phi ,w}$ and every $\rho >0$
the inequalities%
\begin{eqnarray*}
&&\sum_{\mathsf{m}\in \mathbb{Z}}w_{\mathsf{m}}\left( 1+\Phi \left(
\left\vert \mathsf{m}\right\vert \right) \right) ^{k}\Phi \left( \frac{%
\left\vert p_{\mathsf{m}}\right\vert }{\rho }\right) \\
&\leq &\sum_{\mathsf{m}\in \mathbb{Z}}w_{\mathsf{m}}\left( 1+\Phi \left(
\left\vert \mathsf{m}\right\vert \right) \right) ^{k^{\prime }}\Phi \left( 
\frac{\left\vert p_{\mathsf{m}}\right\vert }{\rho }\right) <+\infty
\end{eqnarray*}%
hold true, which proves the desired inclusion.

We proceed by proving the following result, which constitutes in parts a
variant and to some extent a generalization of Proposition 2, Section 10.3
of Chapter X\ in \cite{raoren2}. In what follows we write $l_{\Phi ,w}$ and $%
\mathsf{s}_{\Phi ,w}$ for $l_{k,\Phi ,w}$ and $\mathsf{s}_{k,\Phi ,w}$ when $%
k=0,$ respectively, and $\left\Vert \mathsf{.}\right\Vert _{\Phi ,w}$ for
the corresponding norm:

\bigskip

\textbf{Theorem 3.} \textit{Let }$\Phi $\textit{\ and }$\Psi $\textit{\ be
Orlicz functions. Then the following statements hold:}

\textit{(a) If there exists }$\gamma \in \left( 0,1\right] $ \textit{such
that}%
\begin{equation}
\Phi (t)\leq \Psi (\gamma t)  \label{growthestimate}
\end{equation}%
\textit{for every }$t\in \left[ 0,+\infty \right) $, \textit{and if }$%
k,k^{\prime }\in \mathbb{R}$ \textit{with} 
\begin{equation}
k^{\prime }\geq k\geq 0,  \label{inequality7}
\end{equation}%
\textit{then there exist the continuous embeddings}%
\begin{equation}
l_{k^{\prime },\Psi ,w}\rightarrow l_{k,\Phi ,w}  \label{embedding1}
\end{equation}%
\textit{and}%
\begin{equation}
\mathsf{s}_{k^{\prime },\Psi ,w}\rightarrow \mathsf{s}_{k,\Phi ,w}.
\label{embedding2}
\end{equation}%
\textit{(b) If there exists }$\gamma \in (0,+\infty )$\textit{\ such that}%
\begin{equation}
\Phi (t)\leq \Psi (\gamma t)  \label{inequality8}
\end{equation}%
\textit{for every }$t\in \left[ 0,t_{0}\right] $\textit{\ where} $t_{0}\in
\left( 0,+\infty \right) $, \textit{if }$k\geq 0$ \textit{and} 
\begin{equation}
\inf_{\mathsf{m}\in \mathbb{Z}}w_{\mathsf{m}}>0,  \label{infimumbis}
\end{equation}%
\textit{then there exist the continuous embeddings}%
\begin{equation}
l_{k,\Psi ,w}\rightarrow l_{\Phi ,w}  \label{embedding3}
\end{equation}%
\textit{and}%
\begin{equation}
\mathsf{s}_{k,\Psi ,w}\rightarrow \mathsf{s}_{\Phi ,w}.  \label{embedding4}
\end{equation}

\bigskip

\textbf{Proof. }We first prove that there exists a finite constant $%
c_{\gamma }>0$\ depending solely on $\gamma $\ such that the inequality%
\begin{equation}
\left\Vert \mathsf{p}\right\Vert _{k,\Phi ,w}\leq c_{\gamma }\left\Vert 
\mathsf{p}\right\Vert _{k^{\prime },\Psi ,w}  \label{inequality}
\end{equation}%
holds for every $\mathsf{p}$\ $\in l_{k^{\prime },\Psi ,w}$. According to (%
\ref{luxnorm}), this means that it is sufficient to prove the existence of $%
c_{\gamma }>0$ such that 
\begin{equation}
\sum_{\mathsf{m}\in \mathbb{Z}}w_{\mathsf{m}}\left( 1+\Phi \left( \left\vert 
\mathsf{m}\right\vert \right) \right) ^{k}\Phi \left( \frac{\left\vert p_{%
\mathsf{m}}\right\vert }{c_{\gamma }\left\Vert \mathsf{p}\right\Vert
_{k^{\prime },\Psi ,w}}\right) \leq 1  \label{inequality4}
\end{equation}%
for every non-zero $\mathsf{p}\in l_{k^{\prime },\Psi ,w}$, inequality (\ref%
{inequality}) being obvious if $\mathsf{p=0}$. We recall that for any such $%
\mathsf{p}$ we have%
\begin{equation}
\sum_{\mathsf{m}\in \mathbb{Z}}w_{\mathsf{m}}\left( 1+\Psi \left( \left\vert 
\mathsf{m}\right\vert \right) \right) ^{k^{\prime }}\Psi \left( \frac{%
\left\vert p_{\mathsf{m}}\right\vert }{\left\Vert \mathsf{p}\right\Vert
_{k^{\prime },\Psi ,w}}\right) \leq 1  \label{inequality5}
\end{equation}%
according to (\ref{crucialrelation}). Choosing then any finite constant $%
c\geq \gamma $ and taking (\ref{growthestimate}) and (\ref{inequality7})
into account we have successively 
\begin{eqnarray*}
&&\sum_{\mathsf{m}\in \mathbb{Z}}w_{\mathsf{m}}\left( 1+\Phi \left(
\left\vert \mathsf{m}\right\vert \right) \right) ^{k}\Phi \left( \frac{%
\left\vert p_{\mathsf{m}}\right\vert }{c\left\Vert \mathsf{p}\right\Vert
_{k^{\prime },\Psi ,w}}\right) \\
&\leq &\sum_{\mathsf{m}\in \mathbb{Z}}w_{\mathsf{m}}\left( 1+\Psi \left(
\gamma \left\vert \mathsf{m}\right\vert \right) \right) ^{k}\Psi \left( 
\frac{\gamma \left\vert p_{\mathsf{m}}\right\vert }{c\left\Vert \mathsf{p}%
\right\Vert _{k^{\prime },\Psi ,w}}\right) \\
&\leq &\sum_{\mathsf{m}\in \mathbb{Z}}w_{\mathsf{m}}\left( 1+\Psi \left(
\left\vert \mathsf{m}\right\vert \right) \right) ^{k^{\prime }}\Psi \left( 
\frac{\left\vert p_{\mathsf{m}}\right\vert }{\left\Vert \mathsf{p}%
\right\Vert _{k^{\prime },\Psi ,w}}\right) \leq 1
\end{eqnarray*}%
according to (\ref{inequality5}), which proves (\ref{inequality}) and
thereby embedding (\ref{embedding1}). In order to prove (\ref{embedding2}),
it is sufficient to show that $\mathsf{s}_{k^{\prime },\Psi ,w}\subseteq 
\mathsf{s}_{k,\Phi ,w}$ since $\mathsf{s}_{k^{\prime },\Psi ,w}\subseteq
l_{k^{\prime },\Psi ,w}$. For each $\mathsf{p}\in \mathsf{s}_{k^{\prime
},\Psi ,w}$ and every $\rho >0$ we have%
\begin{eqnarray*}
&&\sum_{\mathsf{m}\in \mathbb{Z}}w_{\mathsf{m}}\left( 1+\Phi \left(
\left\vert \mathsf{m}\right\vert \right) \right) ^{k}\Phi \left( \frac{%
\left\vert p_{\mathsf{m}}\right\vert }{\rho }\right) \\
&\leq &\sum_{\mathsf{m}\in \mathbb{Z}}w_{\mathsf{m}}\left( 1+\Psi \left(
\gamma \left\vert \mathsf{m}\right\vert \right) \right) ^{k}\Psi \left( 
\frac{\gamma \left\vert p_{\mathsf{m}}\right\vert }{\rho }\right) \\
&\leq &\sum_{\mathsf{m}\in \mathbb{Z}}w_{\mathsf{m}}\left( 1+\Psi \left(
\left\vert \mathsf{m}\right\vert \right) \right) ^{k^{\prime }}\Psi \left( 
\frac{\left\vert p_{\mathsf{m}}\right\vert }{\rho }\right) <+\infty ,
\end{eqnarray*}%
which proves the desired inclusion.

We now turn to the proof of Statement (b) by showing that there exists a
finite constant $c_{\gamma ,\delta ,t_{0}}>0$ depending solely on $\gamma $, 
$\delta $ and $t_{0}$ such that the inequality%
\begin{equation}
\left\Vert \mathsf{p}\right\Vert _{\Phi ,w}\leq c_{\gamma ,\delta
,t_{0}}\left\Vert \mathsf{p}\right\Vert _{k,\Psi ,w}  \label{inequality10}
\end{equation}%
holds for every $\mathsf{p}$\ $\in l_{k,\Psi ,w}$, where%
\begin{equation*}
\delta :=\Psi ^{-1}\left( \left( \inf_{\mathsf{m}\in \mathbb{Z}}w_{\mathsf{m}%
}\right) ^{-1}\right)
\end{equation*}%
with $\Psi ^{-1}$ the monotone inverse of $\Psi $. For this it is sufficient
to prove the existence of $c_{\gamma ,\delta ,t_{0}}>0$ implying%
\begin{equation}
\sum_{\mathsf{m}\in \mathbb{Z}}w_{\mathsf{m}}\Phi \left( \frac{\left\vert p_{%
\mathsf{m}}\right\vert }{c_{\gamma ,\delta ,t_{0}}\left\Vert \mathsf{p}%
\right\Vert _{k,\Psi ,w}}\right) \leq 1.  \label{inequality6}
\end{equation}%
From (\ref{infimumbis}) and (\ref{inequality5}) with $k$ instead of $%
k^{\prime }$ we first note that 
\begin{equation*}
\sum_{\mathsf{m}\in \mathbb{Z}}\Psi \left( \frac{\left\vert p_{\mathsf{m}%
}\right\vert }{\left\Vert \mathsf{p}\right\Vert _{k,\Psi ,w}}\right) \leq
\left( \inf_{\mathsf{m}\in \mathbb{Z}}w_{\mathsf{m}}\right) ^{-1}
\end{equation*}%
and \textit{a fortiori}%
\begin{equation*}
\Psi \left( \frac{\left\vert p_{\mathsf{m}}\right\vert }{\left\Vert \mathsf{p%
}\right\Vert _{k,\Psi ,w}}\right) \leq \tilde{c}
\end{equation*}%
for each $\mathsf{m}\in \mathbb{Z}$, where $\tilde{c}=$ $\left( \inf_{%
\mathsf{m}\in \mathbb{Z}}w_{\mathsf{m}}\right) ^{-1}$. Therefore we get 
\begin{equation*}
\frac{\left\vert p_{\mathsf{m}}\right\vert }{c\left\Vert \mathsf{p}%
\right\Vert _{k,\Psi ,w}}\leq \frac{\Psi ^{-1}\left( \tilde{c}\right) }{c}
\end{equation*}%
for every $c>0$. In this manner, requiring 
\begin{equation*}
c\geq \frac{\Psi ^{-1}\left( \tilde{c}\right) }{t_{0}}
\end{equation*}%
we obtain 
\begin{equation*}
\frac{\left\vert p_{\mathsf{m}}\right\vert }{c\left\Vert \mathsf{p}%
\right\Vert _{k,\Psi ,w}}\leq t_{0}
\end{equation*}%
so that we may use (\ref{inequality8}) to get%
\begin{equation*}
\sum_{\mathsf{m}\in \mathbb{Z}}w_{\mathsf{m}}\Phi \left( \frac{\left\vert p_{%
\mathsf{m}}\right\vert }{c\left\Vert \mathsf{p}\right\Vert _{k,\Psi ,w}}%
\right) \leq \sum_{\mathsf{m}\in \mathbb{Z}}w_{\mathsf{m}}\Psi \left( \frac{%
\gamma \left\vert p_{\mathsf{m}}\right\vert }{c\left\Vert \mathsf{p}%
\right\Vert _{k,\Psi ,w}}\right) \text{.}
\end{equation*}%
Consequently, demanding further that $c\geq \gamma $ we have%
\begin{eqnarray*}
&&\sum_{\mathsf{m}\in \mathbb{Z}}w_{\mathsf{m}}\Phi \left( \frac{\left\vert
p_{\mathsf{m}}\right\vert }{c\left\Vert \mathsf{p}\right\Vert _{k,\Psi ,w}}%
\right) \\
&\leq &\sum_{\mathsf{m}\in \mathbb{Z}}w_{\mathsf{m}}\Psi \left( \frac{%
\left\vert p_{\mathsf{m}}\right\vert }{\left\Vert \mathsf{p}\right\Vert
_{k,\Psi ,w}}\right) \leq \sum_{\mathsf{m}\in \mathbb{Z}}w_{\mathsf{m}%
}\left( 1+\Psi \left( \left\vert \mathsf{m}\right\vert \right) \right)
^{k}\Psi \left( \frac{\left\vert p_{\mathsf{m}}\right\vert }{\left\Vert 
\mathsf{p}\right\Vert _{k,\Psi ,w}}\right) \leq 1
\end{eqnarray*}%
according to (\ref{inequality5}), so that the choice%
\begin{equation*}
c\geq \frac{\Psi ^{-1}\left( \tilde{c}\right) }{t_{0}}\vee \gamma
\end{equation*}%
indeed proves (\ref{inequality10}) and thereby embedding (\ref{embedding3}).
The proof of (\ref{embedding4}) follows from a remark similar to that
leading to embedding (\ref{embedding2}). \ \ $\blacksquare $

\bigskip

We are now ready to investigate the existence problem of compact embeddings.
Aside from the decisive role played by the structure of the measure given by
(\ref{simplification}) we show that our analysis rests in an essential way
on the combination of two ingredients, namely, the existence of the Schauder
basis $\left( \mathsf{e}_{\mathsf{m}}\right) $ of Theorem 2 on the one hand,
and the relation%
\begin{equation}
\limsup_{t\rightarrow 0_{+}}\frac{\Phi \left( 2t\right) }{\Phi \left(
t\right) }<+\infty  \label{delta2condition}
\end{equation}%
on the other hand, the latter having been been traditionally referred to as
the $\Delta _{2}$\textit{-condition at zero }in the literature (see, e.g.,
Chapter 4 in \cite{lindenstrausstra}, or \cite{raoren2}). We also prove that
when $k\geq 0$, condition (\ref{delta2condition}) implies that the
Orlicz-Sobolev sequence class coincides with the small and large
Orlicz-Sobolev sequence spaces, respectively. Before stating the precise
result we note the following elementary fact which will play an important
role in its proof, namely, that (\ref{delta2condition}) is equivalent to
having%
\begin{equation*}
\limsup_{t\rightarrow 0_{+}}\frac{\Phi \left( \theta t\right) }{\Phi \left(
t\right) }<+\infty
\end{equation*}%
for every $\theta >0$ in the sense that there exist $t_{\theta }>0$, $%
c_{\theta }>0$ and the bound%
\begin{equation}
\Phi \left( \theta t\right) \leq c_{\theta }\Phi \left( t\right)
\label{bound}
\end{equation}%
for every $t\in \left[ 0,t_{\theta }\right] $. Then we have:

\bigskip

\textbf{Theorem 4. }\textit{The following statements hold:}

\textit{(a) Let }$\Phi $\textit{\ be an Orlicz function satisfying}%
\begin{equation*}
\limsup_{t\rightarrow 0_{+}}\frac{\Phi \left( 2t\right) }{\Phi \left(
t\right) }<+\infty .
\end{equation*}%
\textit{Moreover, let us assume that }$k\geq 0$ \textit{and}%
\begin{equation}
\inf_{\mathsf{m\in }\mathbb{Z}}w_{\mathsf{m}}>0.  \label{infimumter}
\end{equation}%
\textit{Then we have}%
\begin{equation}
\mathsf{s}_{k,\Phi ,w}=\mathsf{c}_{k,\Phi ,w}=l_{k,\Phi ,w},
\label{equalityquarto}
\end{equation}%
\textit{and the Orlicz-Sobolev space in (\ref{equalityquarto}) is a
separable Banach space relative to (\ref{luxnorm}) which consists of all
complex sequences }$\mathsf{p}$\ \textit{such that}%
\begin{equation*}
\sum_{\mathsf{m}\in \mathbb{Z}}w_{\mathsf{m}}\left( 1+\Phi \left( \left\vert 
\mathsf{m}\right\vert \right) \right) ^{k}\Phi \left( \left\vert p_{\mathsf{m%
}}\right\vert \right) <+\infty .
\end{equation*}

\textit{(b) In addition to the above hypotheses, if }$k^{\prime }\in \mathbb{%
R}$ \textit{is such that} \textit{the strict inequality }$k^{\prime }>k$ 
\textit{holds, then embedding (\ref{embeddingbis}) is compact, in which case
we write}%
\begin{equation}
\mathsf{s}_{k^{\prime },\Phi ,w}\hookrightarrow \mathsf{s}_{k,\Phi ,w}.
\label{compactembedding}
\end{equation}

\bigskip

\textbf{Proof. }According to Theorems 1 and 2, it is sufficient to show that 
$l_{k,\Phi ,w}\subseteq \mathsf{s}_{k,\Phi ,w}$ in order to get (\ref%
{equalityquarto}). Let us start with $\mathsf{p}\in l_{k,\Phi ,w}$, so that
there exists at least one $\rho ^{\ast }>0$ ensuring the convergence%
\begin{equation}
\sum_{\mathsf{m}\in \mathbb{Z}}w_{\mathsf{m}}\left( 1+\Phi \left( \left\vert 
\mathsf{m}\right\vert \right) \right) ^{k}\Phi \left( \frac{\left\vert p_{%
\mathsf{m}}\right\vert }{\rho ^{\ast }}\right) <+\infty .
\label{convergenceter}
\end{equation}%
We then have to prove that 
\begin{equation}
\sum_{\mathsf{m}\in \mathbb{Z}}w_{\mathsf{m}}\left( 1+\Phi \left( \left\vert 
\mathsf{m}\right\vert \right) \right) ^{k}\Phi \left( \frac{\left\vert p_{%
\mathsf{m}}\right\vert }{\rho }\right) <+\infty  \label{convergencequarto}
\end{equation}%
for every $\rho >0$. Since $k\geq 0$ and since (\ref{infimumter}) holds we
first infer from (\ref{convergenceter}) that%
\begin{equation*}
\sum_{\mathsf{m}\in \mathbb{Z}}\Phi \left( \frac{\left\vert p_{\mathsf{m}%
}\right\vert }{\rho ^{\ast }}\right) <+\infty ,
\end{equation*}%
which entails the relation%
\begin{equation*}
\lim_{\left\vert \mathsf{m}\right\vert \rightarrow +\infty }\Phi \left( 
\frac{\left\vert p_{\mathsf{m}}\right\vert }{\rho ^{\ast }}\right) =0
\end{equation*}%
and hence 
\begin{equation}
\lim_{\left\vert \mathsf{m}\right\vert \rightarrow +\infty }\frac{\left\vert
p_{\mathsf{m}}\right\vert }{\rho ^{\ast }}=0  \label{limit}
\end{equation}%
from the fact that $\Phi $ is strictly increasing and vanishes only at the
origin. Indeed let $\Phi ^{-1}$ be its monotone inverse, which is
continuous. Then we have%
\begin{equation*}
\Phi ^{-1}\left( \Phi \left( \frac{\left\vert p_{\mathsf{m}}\right\vert }{%
\rho ^{\ast }}\right) \right) =\frac{\left\vert p_{\mathsf{m}}\right\vert }{%
\rho ^{\ast }}\rightarrow \Phi ^{-1}\left( 0\right)
\end{equation*}%
as $\left\vert \mathsf{m}\right\vert \rightarrow +\infty $, but also%
\begin{equation*}
\Phi \left( \Phi ^{-1}\left( 0\right) \right) =0
\end{equation*}%
which implies that $\Phi ^{-1}\left( 0\right) =0$. On the other hand, (\ref%
{bound}) holds for every $\theta >0$ so that there exists $\mathsf{m}%
_{\theta }\in \mathbb{N}^{+}$ implying%
\begin{equation*}
\frac{\left\vert p_{\mathsf{m}}\right\vert }{\rho ^{\ast }}\leq t_{\theta }
\end{equation*}%
for each $\left\vert \mathsf{m}\right\vert \geq \mathsf{m}_{\theta }$
according to (\ref{limit}). Choosing $\theta =\rho ^{-1}\rho ^{\ast }$ we
then have%
\begin{eqnarray*}
&&\sum_{\left\vert \mathsf{m}\right\vert \geq \mathsf{m}_{\theta }}w_{%
\mathsf{m}}\left( 1+\Phi \left( \left\vert \mathsf{m}\right\vert \right)
\right) ^{k}\Phi \left( \frac{\left\vert p_{\mathsf{m}}\right\vert }{\rho }%
\right) \\
&=&\sum_{\left\vert \mathsf{m}\right\vert \geq \mathsf{m}_{\theta }}w_{%
\mathsf{m}}\left( 1+\Phi \left( \left\vert \mathsf{m}\right\vert \right)
\right) ^{k}\Phi \left( \frac{\theta \left\vert p_{\mathsf{m}}\right\vert }{%
\rho ^{\ast }}\right) \\
&\leq &c_{\theta }\sum_{\left\vert \mathsf{m}\right\vert \geq \mathsf{m}%
_{\theta }}w_{\mathsf{m}}\left( 1+\Phi \left( \left\vert \mathsf{m}%
\right\vert \right) \right) ^{k}\Phi \left( \frac{\left\vert p_{\mathsf{m}%
}\right\vert }{\rho ^{\ast }}\right) <+\infty
\end{eqnarray*}%
as a consequence of (\ref{convergenceter}), which proves (\ref%
{convergencequarto}) and thereby the desired inclusion. Separability of (\ref%
{equalityquarto}) then follows from the last statement of Theorem 2.

As for the proof of Statement (b), let $\mathcal{K}$ be a bounded set in $%
\mathsf{s}_{k^{\prime },\Phi ,w}$ and let $\kappa >0$ be the radius of a
ball centered at the origin of $\mathsf{s}_{k^{\prime },\Phi ,w}$ and
containing $\mathcal{K}$, that is, 
\begin{equation*}
\left\Vert \mathsf{p}\right\Vert _{k^{\prime },\Phi ,w}\leq \kappa
\end{equation*}%
for $\mathsf{p}\in \mathcal{K}$. We first show that%
\begin{equation}
\sum_{\mathsf{m}\in \mathbb{Z}}\mu _{k^{\prime },\Phi ,w,\mathsf{m}}\Phi
\left( \frac{\left\vert p_{\mathsf{m}}\right\vert }{\kappa }\right) \leq 1.
\label{inequality11}
\end{equation}%
Indeed, since the preceding statement obviously holds for $\mathsf{p=0}$,
let us assume that $\mathsf{p\neq 0}$ and set momentarily $\mathsf{q:=}%
\kappa ^{-1}\mathsf{p}$ which gives $0<$ $\left\Vert \mathsf{q}\right\Vert
_{k^{\prime },\Phi ,w}\leq 1$. From the convexity of $\Phi $ and the fact
that $\Phi \left( 0\right) =0$ we then get%
\begin{equation*}
\Phi \left( \left\vert q_{\mathsf{m}}\right\vert \right) =\Phi \left(
\left\Vert \mathsf{q}\right\Vert _{k^{\prime },\Phi ,w}\frac{\left\vert q_{%
\mathsf{m}}\right\vert }{\left\Vert \mathsf{q}\right\Vert _{k^{\prime },\Phi
,w}}\right) \leq \left\Vert \mathsf{q}\right\Vert _{k^{\prime },\Phi ,w}\Phi
\left( \frac{\left\vert q_{\mathsf{m}}\right\vert }{\left\Vert \mathsf{q}%
\right\Vert _{k^{\prime },\Phi ,w}}\right)
\end{equation*}%
for every $\mathsf{m}\in \mathbb{Z}$, and therefore%
\begin{eqnarray*}
&&\sum_{\mathsf{m}\in \mathbb{Z}}\mu _{k^{\prime },\Phi ,w,\mathsf{m}}\Phi
\left( \left\vert q_{\mathsf{m}}\right\vert \right) \\
&\leq &\left\Vert \mathsf{q}\right\Vert _{k^{\prime },\Phi ,w}\sum_{\mathsf{m%
}\in \mathbb{Z}}\mu _{k^{\prime },\Phi ,w,\mathsf{m}}\Phi \left( \frac{%
\left\vert q_{\mathsf{m}}\right\vert }{\left\Vert \mathsf{q}\right\Vert
_{k^{\prime },\Phi ,w}}\right) \leq \left\Vert \mathsf{q}\right\Vert
_{k^{\prime },\Phi ,w}
\end{eqnarray*}%
where the last inequality follows from (\ref{crucialrelation}). Switching
back to $\mathsf{p}$ thus gives (\ref{inequality11}). Now from (\ref%
{infimumter}) and the fact that $k^{\prime }\geq 0$ we infer from (\ref%
{inequality11}) the convergence%
\begin{equation*}
\sum_{\mathsf{m}\in \mathbb{Z}}\Phi \left( \frac{\left\vert p_{\mathsf{m}%
}\right\vert }{\kappa }\right) <+\infty ,
\end{equation*}%
which implies%
\begin{equation*}
\lim_{\left\vert \mathsf{m}\right\vert \rightarrow +\infty }\Phi \left( 
\frac{\left\vert p_{\mathsf{m}}\right\vert }{\kappa }\right) =0
\end{equation*}%
and hence 
\begin{equation*}
\lim_{\left\vert \mathsf{m}\right\vert \rightarrow +\infty }\frac{\left\vert
p_{\mathsf{m}}\right\vert }{\kappa }=0
\end{equation*}%
according to the argument already used in the proof of Statement (a).
Therefore, given $\epsilon >0$ and by virtue of (\ref{bound}) with $\theta
=2\epsilon ^{-1}\kappa $, there exist $\mathsf{m}_{1,\epsilon ,\kappa }\in 
\mathbb{N}^{+}$ and constants $c_{\epsilon ,\kappa }>0$, $t_{\epsilon
,\kappa }>0$ such that 
\begin{equation}
\Phi \left( \frac{2\left\vert p_{\mathsf{m}}\right\vert }{\epsilon }\right)
=\Phi \left( \frac{2\kappa \left\vert p_{\mathsf{m}}\right\vert }{\epsilon
\kappa }\right) \leq c_{\epsilon ,\kappa }\Phi \left( \frac{\left\vert p_{%
\mathsf{m}}\right\vert }{\kappa }\right)  \label{inequality12}
\end{equation}%
whenever $\left\vert \mathsf{m}\right\vert \geq \mathsf{m}_{1,\epsilon
,\kappa }$ and $\kappa ^{-1}\left\vert p_{\mathsf{m}}\right\vert \leq
t_{\epsilon ,\kappa }$.

On the other hand we also have%
\begin{equation*}
\lim_{\left\vert \mathsf{m}\right\vert \rightarrow +\infty }\frac{1}{\left(
1+\Phi \left( \left\vert \mathsf{m}\right\vert \right) \right) ^{k^{\prime
}-k}}=0
\end{equation*}%
since $\Phi \left( \left\vert \mathsf{m}\right\vert \right) \rightarrow
+\infty $ as $\left\vert \mathsf{m}\right\vert \rightarrow +\infty $ and $%
k^{\prime }-k>0$. Consequently, there exists $\mathsf{m}_{2,\epsilon ,\kappa
}\in \mathbb{N}^{+}$ such that%
\begin{equation}
\frac{1}{\left( 1+\Phi \left( \left\vert \mathsf{m}\right\vert \right)
\right) ^{k^{\prime }-k}}\leq c_{\epsilon ,\kappa }^{-1}
\label{inequality13}
\end{equation}%
for all $\mathsf{m}\in \mathbb{Z}$ satisfying $\left\vert \mathsf{m}%
\right\vert \geq \mathsf{m}_{2,\epsilon ,\kappa }$, where $c_{\epsilon
,\kappa }$ is the constant appearing in (\ref{inequality12}). Defining then%
\begin{equation}
\mathsf{m}_{\epsilon ,\kappa }:=\mathsf{m}_{1,\epsilon ,\kappa }\vee \mathsf{%
m}_{2,\epsilon ,\kappa }  \label{definition}
\end{equation}%
and using successively (\ref{inequality12}) and (\ref{inequality13}) we
obtain%
\begin{eqnarray}
&&\dsum\limits_{\left\vert \mathsf{m}\right\vert \geq \mathsf{m}_{\epsilon
,\kappa }}\mu _{k,\Phi ,w,\mathsf{m}}\Phi \left( \frac{2\left\vert p_{%
\mathsf{m}}\right\vert }{\epsilon }\right)  \notag \\
&\leq &c_{\epsilon ,\kappa }\dsum\limits_{\left\vert \mathsf{m}\right\vert
\geq \mathsf{m}_{\epsilon ,\kappa }}\frac{1}{\left( 1+\Phi \left( \left\vert 
\mathsf{m}\right\vert \right) \right) ^{k^{\prime }-k}}\mu _{k^{\prime
},\Phi ,w,\mathsf{m}}\Phi \left( \frac{\left\vert p_{\mathsf{m}}\right\vert 
}{\kappa }\right)  \label{estimatebis} \\
&\leq &\dsum\limits_{\left\vert \mathsf{m}\right\vert \geq \mathsf{m}%
_{\epsilon ,\kappa }}\mu _{k^{\prime },\Phi ,w,\mathsf{m}}\Phi \left( \frac{%
\left\vert p_{\mathsf{m}}\right\vert }{\kappa }\right) \leq 1  \notag
\end{eqnarray}%
by virtue of (\ref{inequality11}).

Now every $\mathsf{p=}\left( p_{\mathsf{m}}\right) \in \mathcal{K}$ may be
viewed as an element of $\mathsf{s}_{k,\Phi ,w}$ according to (\ref%
{embeddingbis}). It can therefore be expanded in a unique way as the
norm-convergent series%
\begin{equation*}
\mathsf{p=}\sum_{\mathsf{m}\in \mathbb{Z}}p_{\mathsf{m}}\mathsf{e}_{\mathsf{m%
}}
\end{equation*}%
along the Schauder basis $\left( \mathsf{e}_{\mathsf{m}}\right) $ introduced
in the proof of Theorem 2. Defining 
\begin{equation*}
\mathsf{p}_{\mathsf{m}_{\epsilon ,\kappa }}:\mathsf{=}\sum_{\left\vert 
\mathsf{m}\right\vert <\mathsf{m}_{\epsilon ,\kappa }}p_{\mathsf{m}}\mathsf{e%
}_{\mathsf{m}}
\end{equation*}%
where $\mathsf{m}_{\epsilon ,\kappa }$ is given by (\ref{definition}) we
then have%
\begin{equation*}
\left( \mathsf{p}-\mathsf{p}_{\mathsf{m}_{\epsilon ,\kappa }}\right) _{%
\mathsf{m}}=\left\{ 
\begin{array}{c}
p_{\mathsf{m}}\text{ if }\mathsf{m}\in \mathbb{Z}\cap \left( -\infty \mathsf{%
,}-\mathsf{m}_{\epsilon ,\kappa }\right] \text{ or if }\mathsf{m}\in \text{ }%
\mathbb{Z}\cap \left[ \mathsf{m}_{\epsilon ,\kappa },+\infty \right) \\ 
\\ 
0\text{ otherwise}%
\end{array}%
\right.
\end{equation*}%
for the $\mathsf{m}^{\mathsf{th}}$ component of $\mathsf{p}-\mathsf{p}_{%
\mathsf{m}_{\epsilon ,\kappa }}$, so that%
\begin{equation*}
\left\Vert \mathsf{p}-\mathsf{p}_{\mathsf{m}_{\epsilon ,\kappa }}\right\Vert
_{k,\Phi ,w}=\inf \left\{ \rho >0:\sum_{\left\vert \mathsf{m}\right\vert
\geq \mathsf{m}_{\epsilon ,\kappa }}\mu _{k,\Phi ,w,\mathsf{m}}\Phi \left( 
\frac{\left\vert p_{\mathsf{m}}\right\vert }{\rho }\right) \leq 1\right\}
\end{equation*}%
according to (\ref{luxnorm}). But estimate (\ref{estimatebis}) implies that%
\begin{equation*}
\frac{\epsilon }{2}\in \left\{ \rho >0:\sum_{\left\vert \mathsf{m}%
\right\vert \geq \mathsf{m}_{\epsilon ,\kappa }}\mu _{k,\Phi ,w,\mathsf{m}%
}\Phi \left( \frac{\left\vert p_{\mathsf{m}}\right\vert }{\rho }\right) \leq
1\right\}
\end{equation*}%
and thereby 
\begin{equation*}
\left\Vert \mathsf{p}-\mathsf{p}_{\mathsf{m}_{\epsilon ,\kappa }}\right\Vert
_{k,\Phi ,w}\leq \frac{\epsilon }{2}.
\end{equation*}%
Let us now introduce the finite-dimensional vector space%
\begin{equation*}
\mathcal{S}_{\epsilon ,\kappa }:=\limfunc{span}\left\{ \mathsf{e}_{-\mathsf{m%
}_{\epsilon ,\kappa }},...,\mathsf{e}_{\mathsf{m}_{\epsilon ,\kappa
}}\right\} \subset \mathsf{s}_{k,\Phi ,w}.
\end{equation*}%
Then $\mathsf{p}_{\mathsf{m}_{\epsilon ,\kappa }}\in \mathcal{S}_{\epsilon
,\kappa }$, so that the distance between $\mathsf{p}$ and $\mathcal{S}%
_{\epsilon ,\kappa }$ may be estimated by%
\begin{equation*}
\func{dist}\left( \mathsf{p,}\mathcal{S}_{\epsilon ,\kappa }\right) \leq
\left\Vert \mathsf{p}-\mathsf{p}_{\mathsf{m}_{\epsilon ,\kappa }}\right\Vert
_{k,\Phi ,w}\leq \frac{\epsilon }{2}.
\end{equation*}%
Finally, let $\mathcal{\bar{K}}$ be the closure of $\mathcal{K}$ in $\mathsf{%
s}_{k,\Phi ,w}$ and let $\mathsf{p}\in \mathcal{\bar{K}\setminus K}$. Then
there exists $\mathsf{p}_{\epsilon }\in \mathcal{K}$ such that $\left\Vert 
\mathsf{p}-\mathsf{p}_{\epsilon }\right\Vert _{k,\Phi ,w}\leq \frac{\epsilon 
}{2}$ and therefore%
\begin{equation*}
\func{dist}\left( \mathsf{p,}\mathcal{S}_{\epsilon ,\kappa }\right) \leq 
\frac{\epsilon }{2}+\inf_{\mathsf{q\in }\mathcal{S}_{\epsilon ,\kappa
}}\left\Vert \mathsf{p}_{\epsilon }-\mathsf{q}\right\Vert _{k,\Phi ,w}\leq
\epsilon \text{.}
\end{equation*}%
Thus for every $\mathsf{p}\in \mathcal{\bar{K}}$ we have $\func{dist}\left( 
\mathsf{p,}\mathcal{S}_{\epsilon ,\kappa }\right) \leq \epsilon $, a
statement equivalent to the total boundedness of $\mathcal{\bar{K}}$ (see,
e.g., Proposition 2.1 in \cite{bisgard}). Consequently, $\mathcal{\bar{K}}$
is compact in $\mathsf{s}_{k,\Phi ,w}$. \ \ $\blacksquare $

\bigskip

\textsc{Remark.} One of the key estimates and indeed a turning point in the
above proof is (\ref{estimatebis}), which follows from the $\Delta _{2}$%
-condition at zero thanks to (\ref{inequality12}). But that condition also
implies (\ref{equalityquarto}), so that the question whether the compact
embedding (\ref{compactembedding}) still holds whenever $\mathsf{s}_{k,\Phi
,w}$ is strictly included in $l_{k,\Phi ,w}$ remains open at this time. By
the same token, it is not known whether there exist compact embeddings of
the form $l_{k^{\prime },\Phi ,w}\hookrightarrow l_{k,\Phi ,w}$ whenever $%
k^{\prime }>k$.

\bigskip

Finally, there are many ways to combine the preceding results to generate
new compact embeddings. For instance, the following theorem is an immediate
consequence of the combination of Statements (a) and (b) of Theorems 3 with
Statement (b) of Theorem 4:

\bigskip

\textbf{Theorem 5. }\textit{Let }$\Phi $\textit{\ and }$\Psi $\textit{\ be
Orlicz functions and let us assume that }$\Psi $ \textit{satisfies}%
\begin{equation*}
\limsup_{t\rightarrow 0_{+}}\frac{\Phi \left( 2t\right) }{\Phi \left(
t\right) }<+\infty .
\end{equation*}%
\textit{\ Assume moreover that}%
\begin{equation*}
\inf_{\mathsf{m}\in \mathbb{Z}}w_{\mathsf{m}}>0\text{.}
\end{equation*}%
\textit{Then the following statements hold:}

\textit{(a) If the inequality}%
\begin{equation*}
\Phi (t)\leq \Psi (\gamma t)
\end{equation*}%
\textit{holds} \textit{for some }$\gamma \in \left( 0,1\right] $ \textit{and
every }$t\in \left[ 0,+\infty \right) $,\textit{\ and if }$k^{\prime \prime
}>k^{\prime }\geq k\geq 0$, \textit{then we have the embeddings}%
\begin{equation*}
\mathsf{s}_{k^{\prime \prime },\Psi ,w}\hookrightarrow \mathsf{s}_{k^{\prime
},\Psi ,w}\rightarrow \mathsf{s}_{k,\Phi ,w}
\end{equation*}%
\textit{where the first one is compact and the second one continuous. In
particular, the embedding}%
\begin{equation*}
\mathsf{s}_{k^{\prime \prime },\Psi ,w}\hookrightarrow \mathsf{s}_{k,\Phi ,w}
\end{equation*}%
\textit{is compact.}

\textit{(b) If the inequality}%
\begin{equation*}
\Phi (t)\leq \Psi (\gamma t)
\end{equation*}%
\textit{holds} \textit{for some }$\gamma \in \left( 0,+\infty \right) $%
\textit{\ and every }$t\in \left[ 0,t_{0}\right] $\textit{\ with }$t_{0}\in
\left( 0,+\infty \right) $,\textit{\ and if }$k^{\prime }>k\geq 0$\textit{,
then we have the embeddings}%
\begin{equation*}
\mathsf{s}_{k^{\prime },\Psi ,w}\hookrightarrow \mathsf{s}_{k,\Psi
,w}\rightarrow \mathsf{s}_{\Phi ,w}
\end{equation*}%
\textit{where the first one is compact and the second one continuous. In
particular, the embedding}%
\begin{equation*}
\mathsf{s}_{k^{\prime },\Psi ,w}\hookrightarrow \mathsf{s}_{\Phi ,w}
\end{equation*}%
\textit{is compact.}

\bigskip

\textsc{Remark.} It is worth making one last observation regarding
compactness. It was indeed proven in \cite{vuillermot} that the spaces $h_{%
\mathbb{C},w}^{k,s}$ enjoy Pitt's property, namely, the fact that if $%
s>s^{\prime }\geq 1$ then every linear bounded operator from $h_{\mathbb{C}%
,w}^{k,s}$ into $h_{\mathbb{C},w}^{k,s^{\prime }}$ is compact, thereby
generalizing Pitt's original result which was about the usual sequence
spaces $l^{s}$ and $l^{s^{\prime }}$ corresponding to $k=0$ and $w_{\mathsf{m%
}}=1$ for every $\mathsf{m}\in \mathbb{Z}$. Pitt's original proof appeared
in \cite{pitt} and was later considerably simplified in \cite{delpech} and 
\cite{fabianzizler}. A glance at the proof of Theorem 3 in \cite{vuillermot}
shows, however, that the arguments set forth there are not applicable to the
more general spaces introduced in this article. An interesting open problem
is therefore to see under what conditions on the Orlicz functions $\Phi $
and $\Psi $ Pitt's property might hold for $\mathsf{s}_{k,\Phi ,w}$ and $%
\mathsf{s}_{k,\Psi ,w}$, perhaps in the spirit of the remark following the
proof of Theorem 4.b.3 in Chapter 4 of \cite{lindenstrausstra} corresponding
to the case $k=0$.

\bigskip

We devote the next section to illustrating some of the preceding statements.

\section{Some examples}

1. Let us consider the functions $\Phi \left( t\right) =t^{r}$ and $\Psi
(t)=t^{s}$ where $r,s\in \left[ 1,+\infty \right) $, which both satisfy the $%
\Delta _{2}$-condition at zero. Let $w=$ $\left( w_{\mathsf{m}}\right) $ be
a sequence of weights satisfying (\ref{infimumter}) and let $k\geq 0$. Then%
\begin{equation*}
\mathsf{s}_{k,\Phi ,w}=l_{k,\Phi ,w}=h_{\mathbb{C},w}^{k,r}
\end{equation*}%
and%
\begin{equation*}
\mathsf{s}_{k,\Psi ,w}=l_{k,\Psi ,w}=h_{\mathbb{C},w}^{k,s},
\end{equation*}%
where $h_{\mathbb{C},w}^{k,r}$ and $h_{\mathbb{C},w}^{k,s}$ are the Banach
spaces of all complex sequences satisfying (\ref{norm}) with $r$ and $s,$
respectively. Let us assume that $r>s$ and for a given $\gamma >0$ let us
set $t_{0}=\gamma ^{\frac{s}{r-s}}$. Then we have%
\begin{equation*}
\Phi (t)\leq \Psi (\gamma t)
\end{equation*}%
for every $t\in \left[ 0,t_{0}\right] $. Thus with $k^{\prime }>k\geq 0$ and%
\begin{equation*}
\inf_{\mathsf{m\in }\mathbb{Z}}w_{\mathsf{m}}>0
\end{equation*}%
we have the embeddings%
\begin{equation*}
h_{\mathbb{C},w}^{k^{\prime },s}\hookrightarrow h_{\mathbb{C}%
,w}^{k,s}\rightarrow h_{\mathbb{C},w}^{0,r}
\end{equation*}%
where the first one is compact and the second one continuous, as a direct
application of Statement (b) in Theorem 5. In particular, the embedding%
\begin{equation*}
h_{\mathbb{C},w}^{k^{\prime },s}\hookrightarrow h_{\mathbb{C},w}^{0,r}
\end{equation*}%
is compact. In this way we retrieve some of the results proved in \cite%
{vuillermot}.

\bigskip

2. This example and the next one are a direct application of Theorem 4. Let
us consider the Orlicz function\textit{\ }%
\begin{equation}
\Phi \left( t\right) =\exp \left[ t^{2}\right] -1  \label{orliczfunction}
\end{equation}%
which satisfies the $\Delta _{2}$-condition at zero. If $k\geq 0$ and%
\begin{equation*}
\inf_{\mathsf{m\in }\mathbb{Z}}w_{\mathsf{m}}>0,
\end{equation*}%
then the Orlicz-Sobolev space $\mathsf{s}_{k,\Phi ,w}=$\ $l_{k,\Phi ,w}$ is
a separable Banach space relative to (\ref{luxnorm}) which consists of all
complex sequences $\mathsf{p}$\ such that 
\begin{equation*}
\sum_{\mathsf{m}\in \mathbb{Z}}w_{\mathsf{m}}\exp \left[ k\mathsf{m}^{2}%
\right] \left( \exp \left[ \left\vert p_{\mathsf{m}}\right\vert ^{2}\right]
-1\right) <+\infty .
\end{equation*}%
Moreover, if $k^{\prime }\in \mathbb{R}$ is such that the strict inequality $%
k^{\prime }>k$ holds, then the embedding%
\begin{equation*}
\mathsf{s}_{k^{\prime },\Phi ,w}\hookrightarrow \mathsf{s}_{k,\Phi ,w}
\end{equation*}%
is compact.

\bigskip

3.\ Under the same hypotheses, identical conclusions hold for the Orlicz
function\textit{\ }%
\begin{equation}
\Phi \left( t\right) =\exp \left[ t\right] -t-1,  \label{orliczfunctionbis}
\end{equation}%
in which case the Orlicz-Sobolev space $\mathsf{s}_{k,\Phi ,w}=l_{k,\Phi ,w}$
consists of all complex sequences $\mathsf{p}$\ such that 
\begin{equation*}
\sum_{\mathsf{m}\in \mathbb{Z}}w_{\mathsf{m}}\left( \exp \left[ \left\vert 
\mathsf{m}\right\vert \right] -\left\vert \mathsf{m}\right\vert \right)
^{k}\left( \exp \left[ \left\vert p_{\mathsf{m}}\right\vert \right]
-\left\vert p_{\mathsf{m}}\right\vert -1\right) <+\infty .
\end{equation*}

\bigskip

4. Let $\Phi $ be an Orlicz function satisfying the $\Delta _{2}$-condition
at zero and let $\Psi :\left[ 0,+\infty \right) \rightarrow \left[ 0,+\infty
\right) $ be the function defined by%
\begin{equation}
\Psi (t):=\exp \left[ \Phi (t)\right] -1.  \label{orliczfunctionter}
\end{equation}%
It is then easily verified that $\Psi $ is also an Orlicz function
satisfying the $\Delta _{2}$-condition at zero, and that%
\begin{equation*}
\Phi (t)\leq \Psi (\gamma t)
\end{equation*}%
for each $\gamma \in \left[ 1,+\infty \right) $ and every $t\in \left[
0,t_{0}\right] $ for an arbitrarily chosen $t_{0}\in \left( 0,+\infty
\right) .$ Thus with $k^{\prime }>k\geq 0$ and%
\begin{equation*}
\inf_{\mathsf{m\in }\mathbb{Z}}w_{\mathsf{m}}>0
\end{equation*}%
we have%
\begin{equation*}
\mathsf{s}_{k,\Phi ,w}=\ l_{k,\Phi ,w},
\end{equation*}%
\begin{equation*}
\mathsf{s}_{k,\Psi ,w}=\ l_{k,\Psi ,w},
\end{equation*}%
and the conclusions of Statement (b) of Theorem 5 hold.\textsf{\ }In
particular, the embedding%
\begin{equation*}
\mathsf{s}_{k^{\prime },\Psi ,w}\hookrightarrow \mathsf{s}_{\Phi ,w}
\end{equation*}%
is compact. For instance, if $\Phi \left( t\right) =t^{r}$ and $\Psi
(t)=\exp \left[ t^{r}\right] -1$ with $r\in \left[ 1,+\infty \right) $ then
the embedding%
\begin{equation*}
\mathsf{s}_{k^{\prime },\Psi ,w}\hookrightarrow h_{\mathbb{C},w}^{0,r}
\end{equation*}%
is compact.

\bigskip

\textsc{Remark. }The properties of the Orlicz-Sobolev sequence spaces
described in the last three examples are essentially different from those of
the corresponding spaces defined with respect to Lebesgue measure in regions
of Euclidean space. Indeed, in the latter case the exponential growth of (%
\ref{orliczfunction}), (\ref{orliczfunctionbis}) and (\ref{orliczfunctionter}%
) implies that the corresponding Orlicz-Sobolev classes \textit{do} \textit{%
not} carry a linear structure, thereby making the small and the large
Orlicz-Sobolev spaces quite different objects (see, e.g., Chapter VIII in 
\cite{adamsfournier}, or \cite{krasnorut}). More generally, the dichotomy
between certain results obtained for spaces defined with respect to
diffusive measures such as Lebesgue's and those valid for spaces generated
by counting or discrete measures has been amply discussed in \cite{raoren1}
and \cite{raoren2}. We refer the reader to these references for more details.

\bigskip

\textbf{Statements and declarations}

\bigskip

\textsc{Declaration of competing interest:} The author declares that he has
no known competing financial interests or personal relationships that could
have appeared to influence the work reported in this paper.

\bigskip

\textsc{Conflict of interest:} The author states that there is no conflict
of interest regarding the content of this paper.

\bigskip

\textsc{Data availability statement:} All the data supporting the results
stated in this article are available in the bibliographical references liste
below.

\bigskip

\textsc{Acknowledgements:}\textbf{\ }The author would like to thank the Funda%
\c{c}\~{a}o para a Ci\^{e}ncia e a Tecnologia (FCT) of the Portuguese
Government for its financial support under grant UIDB/04561/2020.

\end{document}